\documentclass[12pt]{article}
\usepackage{amsmath}
\usepackage{amssymb}
\usepackage{theorem}
\newtheorem{theorem}{Theorem}
\newtheorem{cor}[theorem]{Corollary}
\newtheorem{prop}[theorem]{Proposition}
\newtheorem{lemma}[theorem]{Lemma}
\newenvironment{proof}{\textbf{Proof.}}{\hfill $\square$ \\}
\newcommand{\vol}{\text{vol}^2\hspace{-.5mm}}
\newtheorem{remark}[theorem]{Remark}
\newcommand{\lf}{L\!f}
\newcommand{\floor}[1]{\lfloor #1 \rfloor}
\newcommand{\inprod}[2]{\langle #1 \, ,#2 \rangle}
\newcommand{\inprodm}[2]{\langle #1 \, ,#2 \rangle_M}
\newcommand{\foral}{\quad \text{for all} \quad}
\newenvironment{defin}{\textbf{Definition.}}{}
\newcommand{\sph}{S^{n-1}}
\newcommand{\li}[2]{||\hspace{.5mm}l_{#1}(#2)\hspace{-.35mm}||^{\hspace{.1mm}2}}
\begin{document}
{\title{Convexity Properties of The Cone of Nonnegative
Polynomials}}
\author{Grigoriy Blekherman}
\date{} \maketitle \numberwithin{theorem}{section}
\begin{abstract}
\noindent We study metric properties of the cone of homogeneous
non-negative multivariate polynomials and the cone of sums of
powers of linear forms, and the relationship between the two
cones. We compute the maximum volume ellipsoid of the natural base
of the cone of non-negative polynomials and the minimum volume
ellipsoid of the natural base of the cone of powers of linear
forms and compute the coefficients of symmetry of the bases. The
multiplication by $(x_1^2 + \ldots + x_n^2)^m$ induces an
isometric embedding of the space of polynomials of degree $2k$
into the space of polynomials of degree $2(k+m)$, which allows us
to compare the cone of non-negative polynomials of degree $2k$ and
the cone of sums of $2(k+m)$-powers of linear forms. We estimate
the volume ratio of the bases of the two cones and the rate at
which it approaches 1 as $m$ grows.

\end{abstract}
\section{Introduction and Results}
\subsection{Introduction}
\hspace{.44cm} Let $P_{n,d}$ denote the vector space of real
homogeneous polynomials (forms) of degree $d$ in $n$ real
variables. For even $d=2k$ there are three interesting closed
convex cones in $P_{n,2k}$:
\\
\noindent The cone of nonnegative polynomials, $C(=C_{n,2k})$
\begin{equation*}
C=\bigl{\{}f \in P_{n,2k} \mid f(x) \geq 0 \quad \text{for all}
\quad x \in \mathbb{R}^n \bigr{\}}.
\end{equation*}
The cone of sums of squares, $Sq(=Sq_{n,2k})$
\begin{equation*}
Sq=\biggl{\{} f \in P_{n,2k} \mathrel{\bigg{\arrowvert}} f=\sum_i
f_i^2 \quad \text{for some} \quad f_i \in P_{n,k} \biggr{\}}.
\end{equation*}
The cone of sums of $2k$-th powers of linear forms,
$\lf(=\lf_{n,2k})$
\begin{equation*}
\lf=\biggl{\{} f \in P_{n,2k} \mathrel{\bigg{\arrowvert}} f=\sum_i
l_i^{2k} \quad \text{for some linear forms} \quad l_i \in
P_{n,1}\bigg{\}}.
\end{equation*}
\\ \indent The study of algebraic properties of these cones goes
back to Hilbert, who described explicitly all the cases when
$C_{n,2k}=Sq_{n,2k}$,\cite{hilbert}. Hilbert's 17th problem,
solved in affirmative by Artin and Schreier in the 1920's, asked
whether every nonnegative polynomial is a sum of squares of
rational functions \cite{realal}. Constructive aspects of
Hilbert's problem still draw attention today
\cite{realal},\cite{rez1}. For a discussion of some algebraic
properties of the cone of sums of powers of linear forms we refer
to \cite{rez2}.
\\\indent To our knowledge, however, these cones have not been
studied as general convex objects, possessing invariants based on
convexity. In this paper we look at some convex properties of
these cones.
\\ \indent Let $M(=M_{n,2k})$ denote the hyperplane of all
forms in $P_{n,2k}$ with integral 1 on the unit sphere $\sph$:
\begin{equation*}
M=\biggr{\{} f \in P_{n,2k} \mathrel{\bigg{\arrowvert}}
\int_{\sph}f \, d\sigma=1 \biggl{\}},
\end{equation*}
where $\sigma$ denotes the rotation invariant probability measure
on $\sph$.

We define compact convex bodies $\overline{C}$,
$\overline{Sq}$ and
$\overline{\lf}$
by intersecting the respective cones with $M$:
\begin{equation*}
\overline{C}=C \cap M, \quad \overline{Sq}=Sq \cap M, \quad
\text{and} \quad \overline{\lf}=\lf \cap M.
\end{equation*}
The compact convex bodies $\overline{C}$, $\overline{Sq}$ and
$\overline{\lf}$ are natural bases of the respective cones and
they have full dimension in $M$. Their naturality becomes apparent
if we consider the following action of the special orthogonal
group $SO(n)$ on $P_{n,d}$:
\begin{equation*}
A \in SO(n) \quad \text{sends} \quad  f(x) \in P_{n,d} \quad
\text{to} \quad Af=f(A^{-1}x).
\end{equation*}
All three cones $C$, $Sq$ and $\lf$ are fixed by the action of
$SO(n)$, and $M$ is the only hyperplane in $P_{n,2k}$ fixed by
this action. Therefore $\overline{C}$, $\overline{Sq}$ and
$\overline{\lf}$ are also fixed by the action of $SO(n)$, and they
are the only hyperplane sections of their respective cones with
this property. This action of $SO(n)$ naturally gives a
homomorphism
\begin{equation*}
\phi_{n,d}:SO(n) \rightarrow GL(P_{n,d}),
\end{equation*}
and therefore we have a representation of $SO(n)$ on $P_{n,d}$.
There is a natural inner product on $P_{n,d}$:
\begin{equation*}
\inprod{f}{g}=\int_{S^{n-1}} fg \, d\sigma.
\end{equation*}
The metric induced by this inner product makes $\phi_{n,d}$ an
orthogonal representation, since the inner product is invariant
under the action of $SO(n)$.
\\ \indent Let $K(2m)$ be the cone in $P_{n,2k}$ of forms whose restrictions
to the sphere are linear combinations of $2m$-th powers of linear
forms on $\sph$. Equivalently $K(2m)$ is the cone of forms in
$P_{n,2k}$ that multiplied by $(x_1^2+\ldots +x_n^2)^{m-k}$ become
sums of powers of linear forms
\begin{equation*}
K(2m)=\biggl{\{} f \in P_{n,2k} \mathrel{\bigg{\arrowvert}}
(x_1^2+\ldots+x_n^2)^{m-k}f \in \lf_{n,2m} \biggl{\}}.
\end{equation*}
We define $\overline{K}(2m)$ by intersecting $K$ with the
hyperplane of forms of integral 1 on $\sph$.
\\ \indent From general convexity we know that every compact convex body $K$
contains a unique ellipsoid of maximum volume, known as John's
ellipsoid of $K$. Also, $K$ is contained in a unique ellipsoid of
minimum volume, known as the Loewner ellipsoid of $K$,
\cite{ball}.
\\ \indent A crude, yet interesting, measure of
symmetry of $K$ is its coefficient of symmetry about a point $v$
in the interior of $K$. The coefficient of symmetry of $K$ about
$v$ is defined as the largest $\alpha \in \mathbb{R}$ such that
\begin{equation*}
-\alpha(K-v) \quad  \text{ is contained in} \quad K-v .
\end{equation*}
We will compute coefficients of symmetry of $\overline{C}$ and
$\overline{\lf}$ with respect to $v=(x_1^2+\ldots+x_n^2)^k$.
\subsection{Convexity Results}
We prove the following properties:
\begin{enumerate}
\item Let $\alpha=\text{dim}P_{n,2k}-1$. Then
\begin{equation*}
\left( \frac{\text{vol}\,\overline{K}(2m)}{\text{vol}\,
\overline{C}} \right)^{1/\alpha} \geq
\frac{m!\Gamma(\frac{2m+n}{2})} {(m-k)!\Gamma(\frac{2m+2k+n}{2})}.
\end{equation*}
It follows that if for an $\epsilon >0$ we let
$m=(2k^2+kn)/{\epsilon}$ then
\begin{equation*}
\left(\frac{\text{vol}\,\overline{K}(2m)}{\text{vol}\,
\overline{C}} \right)^{1/\alpha} \geq 1-\epsilon.
\end{equation*}
(cf Theorem \ref{volumeratio}). Thus the volume ratio approaches 1
as $m$ goes to infinity. Therefore all strictly positive
polynomals lie in some $K(2m)$. (cf \cite{rez1})
\item We show
that, in the above metric, John's ellipsoid of
$\overline{C}_{n,2k}$ is a ball centered at $(x_1^2+ \ldots
+x_n^2)^k$ of radius
\begin{equation*}
\frac{1}{\sqrt{\binom{n+2k-1}{2k}-1}}.
\end{equation*}
(cf Theorem \ref{main}). \item \label{lf}We explicitly compute the
Loewner Ellipsoid of $\overline{\lf}_{n,2k}$. (cf Theorem
\ref{mainpowers}).
\item We calculate the coefficient of symmetry
of $\overline{C}_{n,2k}$ and $\overline{\lf}_{n,2k}$ with respect
to $(x_1^2+ \ldots +x_n^2)^k$, which in both cases turns out to be
\begin{equation*}
\frac{1}{\binom{n+k-1}{k}-1}.
\end{equation*}
(cf Theorem \ref{nonnegcoeffsym} and Theorem \ref{mainpowers}).
\item Combining (2) and (3) we show that $\overline{\lf}_{n,2k}$
contains a ball of radius
\begin{equation*}
\frac{k!\Gamma(k+\frac{n}{2})}{\Gamma(2k+\frac{n}{2})\sqrt{\binom{n+k-1}{k}-1}},
\end{equation*}
centered at $(x_1^2+\ldots +x_n^2)^k$.  (cf Corollary
\ref{powerball}). \item A crucial tool for the above calculations
is computation of the Loewner ellipsoid of a convex hull of the
orbit of an arbitrary point under a continuous group action of a
compact group. We apply this to the case of $SO(n)$. (cf Theorem
\ref{orbits}).
\end{enumerate}

Note that the invariants computed in (2)-(4) are independent of
Euclidean structure on $P_{n,d}$, e.g. the maximal volume
ellipsoid is unique and is the same regardless of the choice of an
inner product, although it will not always be a ball.
\\ \indent In many cases we reduce our calculations to polynomials symmetric
with respect to an axis. These are the polynomials fixed by
$J(n,v)$, where $J(n,v)$ is the subgroup of $SO(n)$ consisting of
orthogonal transformations that fix a particular vector $v \in
\mathbb{R}^n$:
\begin{equation*}
J(n,v)=\{A \in SO(n) \mid Av=v \quad \text{for some fixed} \quad v
\in \mathbb{R}^n\}.
\end{equation*}
We show that every nonnegative polynomial symmetric with respect
to an axis is a sum of squares, which proves to be quite useful,
and we think interesting in itself.
\subsection{Integral Inequalities}
\hspace{.44cm} A byproduct of our work is a number of integral
inequalities for homogeneous polynomials on $\sph$. We use the
usual notation for $L^p$ and $L^{\infty}$ norms:
\begin{equation*}
||f||_p=\left(\int_{\sph} |f|^p d\sigma \right)^{\frac{1}{p}}
 \quad \text{and} \quad ||f||_{\infty}=\max_{x \in S^{n-1}} |f(x)|.
\end{equation*}
We list some of our results roughly in order of appearance in the
paper:
\begin{enumerate}



\item \label{nifty} For nonnegative $f \in P_{n,2k}$,
\begin{equation*}
||f||_{\infty} \leq \binom{n+k-1}{k}||f||_1.
\end{equation*}
(cf Theorem \ref{maxi3})
\item \label{nifty2} Equivalently to
\eqref{nifty}, let $M_f$ denote the maximum of $f$ on $\sph$ and
$m_f$ denote the minimum of $f$ on $\sph$. Also let
$\alpha=\frac{1}{\binom{n+k-1}{k}}$. Then for all $f \in P_{n,2k}$

\begin{equation*}
\alpha M_f + (1-\alpha) m_f \leq \int_{\sph} f \, d\sigma \leq
(1-\alpha)M_f+\alpha m_f.
\end{equation*}
(cf Corollary \ref{ineqcor})
\item \label{n2} For $f$ as in
\eqref{nifty},
\begin{equation*}
||f||_2 \leq \sqrt{\binom{n+k-1}{k}}||f||_1.
\end{equation*}
(cf Corollary \ref{mainint})

\item \label{niftygen} From \eqref{nifty} we easily derive that
for $f \in P_{n,k}$
\begin{equation*}
||f||_{\infty} \leq {\binom{n+kl-1}{kl}}^{\frac{1}{2l}}||f||_{2l}
\, ,
\end{equation*}
for all positive integers $l$.
(cf Corollary \ref{ineqcor1})
\end{enumerate}
Estimates $(1)-(3)$ above are sharp and we also provide all
extreme polynomials for them. For a different proof of
\eqref{niftygen} by Barvinok and a discussion of applications see
\cite{barv}. Sogge in \cite{sogge}, and Duoandikoetxea in
\cite{duo} derive some related interesting inequalities.
\\ \indent The rest of the article is structured as follows:
Section 2 contains the known results necessary for the rest of the
paper. In Section 3 we compute the Loewner ellipsoid of an orbit
of a point under the action of a compact group. In Section 4 we
prove some results about polarity in the space of forms with
respect to our inner product $\inprod{}{}$. In Section 5 we
compute John's ellipsoid for the cone of nonnegative polynomials.
In Section 6 we compute the coefficient of symmetry of the cone of
nonnegative polynomials. Section 7 is devoted to the cone of sums
of powers of linear forms. We derive the equation of its Loewner
ellipsoid and compute its coefficient of symmetry, and we show the
volume ratio result.
\section{Preliminaries}
\subsection{Representation of $SO(n)$ in $P_{n,d}$}
\hspace{.44cm} There is a natural action of $SO(n)$ on $P_{n,d}$
which sends $f(x)$ to $f(A^{-1}x)$ for $A \in SO(n)$. We will
denote the action of $A \in SO(n)$ on $f$ by $Af$. Note that this
leads naturally to a representation $\phi_{n,d}$ of $SO(n)$. We
introduce an inner product on $P_{n,d}$:
\begin{equation*}
\inprod {f}{g} =\int_{S^{n-1}} f(x)g(x) \, d\sigma
\end{equation*}
where $\sigma$ denotes the rotation invariant probability measure
on $S^{n-1}$. Under our inner product the norm of $f$ coincides
with the usual $L^2$ norm and we will often use $||f||$ instead of
$||f||_2$. The metric induced by the inner product makes
$\phi_{n,d}$ into an orthogonal representation as
\begin{equation*}
\inprod{Af}{Ag}=\int_{S^{n-1}} f(A^{-1}x)g(A^{-1}x) d\sigma=
\int_{S^{n-1}} f(x)g(x) d\sigma=\inprod{f}{g},
\end{equation*}
by rotational invariance of $\sigma$.
\\ \indent We use $\Delta$
to denote the Laplace differential operator:
\begin{equation*}
\Delta=\frac{\partial^2}{\partial
x_1^2}+\cdots+\frac{\partial^2}{\partial x_n^2}.
\end{equation*}
\begin{defin}
If $f \in P_{n,d}$ and
\begin{equation*}
\Delta f = 0,
\end{equation*}
then $f$ is called a
homogeneous harmonic.
\end{defin}
\\\indent The restriction of a homogeneous harmonic to the sphere
$\sph$ is called a spherical harmonic. By linearity of $\Delta$,
homogeneous harmonics form a vector subspace of $P_{n,d}$, which
we denote by $H_{n,d}$:
\begin{equation*}
H_{n,d}=\{f \in P_{n,d} \mid \Delta f=0\}.
\end{equation*}
Let
\begin{equation*}
r(x)=(x_1^2+ \ldots +x_n^2)^{1/2}.
\end{equation*}
The inclusion $i: H_{n,d-2l} \to P_{n,d}$ given by
\begin{equation*}
i(f)=r^{2l}f
\end{equation*}
is an isometry since $i(f)=r^{2l}f$ is the same function as $f$ on
the sphere $\sph$. We denote the image subspace of $P_{n,d}$ by
$H^*_{n,d-2l}$:
\begin{equation*}
H^*_{n,d-2l}=\{f \in P_{n,d} \mid f=r^{2l}g \ \text{for some} \ g
\in H_{n,d-2l}\}
\end{equation*}
\\We need some facts about the representations $\phi_{n,d}$, see \cite{mull} and \cite{vilenkin}.
\begin{theorem}
\label{represent} $H_{n,d}$ is an irreducible $SO(n)$-module, and,
therefore, $H^*_{n,d}$ is an irreducible submodule of $P_{n,d}$.
Furthermore, $P_{n,d}$ splits into irreducible submodules as
follows:
\begin{equation*}
P_{n,d}=\bigoplus_{i=0}^{\lfloor d/2 \rfloor}
r^{2i}H_{n,d-2i}=\bigoplus_{i=0}^{\lfloor d/2 \rfloor}
H^*_{n,d-2i}.
\end{equation*}
Let $D(n,d)$ be the dimension of $P(n,d)$ and let $N(n,d)$ be the
dimension of $H_{n,d}$. Then
\begin{equation*}
D(n,d)=\binom{n+d-1}{d} \quad \text{and} \quad
N(n,d)=\frac{(2d+n-2)(d+n-3)!}{d!(n-2)!}.
\end{equation*}
\end{theorem}
\begin{remark}
The restriction of $f \in P_{n,d}$ to the sphere $\sph$ can be
uniquely written as a sum of spherical harmonics of degrees having
the same parity as $d$.
\end{remark}
\begin{defin}
Let $J(n,v)$ denote the subgroup of $SO(n)$ that keeps a
particular $v \in \sph$ fixed:
\begin{equation*}
J(n,v)=\{A \in SO(n) \mid Av=v\}.
\end{equation*}
\end{defin}
We denote the standard basis of $\mathbb{R}^n$ by $e_1 \ldots
e_n$. We will use the following theorem on restricting
$\phi_{n,d}$ to $J(n,v)$:
\begin{theorem}
\label{superbasic} There exists unique polynomial $L^v_{n,d}(x)$
with the following properties:
\begin{enumerate}
\item $L^v_{n,d}(x) \in H_{n,d}$, \item
$L^v_{n,d}(Ax)=L^v_{n,d}(x)$ for all $A \in J(n,v)$, \item
$L^v_{n,d}(v)=1$.
\end{enumerate}
We will call $L^v_{n,d}(x)$ the Legendre harmonic with axis $v$.
(also called axial and zonal polynomial).
\end{theorem}
We will denote $L_{n,d}^{e_n}$ simply by $L_{n,d}$. We now state
some facts about Legendre harmonics that will be used later on:
\begin{theorem} \hspace{1pt}
\\
\label{facts}
\begin{enumerate}
\item The norm of the Legendre harmonic is given by:
\begin{equation*}
||L^v_{n,d}||^2=\int_{S^{n-1}} (L^v_{n,d})^2 \ d\sigma=
\frac{1}{N(n,d)}
\end{equation*}
\item
\begin{equation*}
||L^v_{n,d}||_{\infty}=1.
\end{equation*}
The maximum absolute value of $L^v_{n,d}$ is achieved only at $v$,
$-v$, and
\begin{equation*}
L^v_{n,2k}(v)=L^v_{n,2k}(-v)=1 \quad \text{while} \quad
L^v_{n,2k-1}(v)=-L^v_{n,2k-1}(-v)=1.
\end{equation*}
\end{enumerate}
\end{theorem}

\indent Since $L_{n,d}$ is fixed by $J(n,e_n)$, by applying
rotations of $\sph$ fixing $e_n$, we see that $L_{n,d}$ is
constant on slices of the sphere with hyperplanes $T_a$
perpendicular to $e_n$:
\begin{equation*}
T_{a}=\big{\{}\zeta \in \mathbb{R}^n \mid \inprod{\zeta}{e_n}=a
\big{\}}, \quad -1 \leq a \leq 1.
\end{equation*}
Hence the Legendre harmonics $L_{n,d}$ restricted to the sphere
$\sph$ are functions of essentially only one variable, namely, the
last coordinate. Therefore we can define a polynomial in $t$,
which we will denote $Q_{n,d}(t)$, such that
\begin{equation*}
L_{n,d}(\xi)=Q_{n,d}(\langle \xi, e_n \rangle) \quad \text{for
all} \quad \xi \in S^{n-1}.
\end{equation*}
The family of polynomials $Q_{n,d}(t)$ are known as the Legendre
polynomials and are special cases of ultraspherical (or
Gegenbauer) polynomials. For many identities satisfied by these
polynomials see \cite{szeg} and \cite{vilenkin}.
\subsection{Loewner and John Ellipsoids}
\hspace{.44cm} Let $K$ be a convex body in a finite dimensional
real vector space $V$. There exists a unique ellipsoid of maximal
volume contained in $K$, known as John's ellipsoid of $K$; we will
denote it by $D_K$. Moreover, there is a criterion for determining
whether a given ellipsoid $E$ contained in $K$ is John's ellipsoid
of $K$ based solely on the points in the intersection of
boundaries $\,\partial E \cap \partial K$.
\\\indent Recall
that a non-singular linear transformation does not affect ratios
of volumes. Therefore, after translating the center of $D_K$ to
the origin and then applying a linear transformation $A \in
GL(V)$, we know that John's Ellipsoid of $A(K)$ is the unit ball
$B^n$. Therefore we will assume that John's ellipsoid of $K$ is a
ball and we state the theorem for this case:
\begin{theorem}
Each convex body $K$ contains unique ellipsoid of maximal volume.
This ellipsoid is $B^n$ if and only if the following conditions
are satisfied: $B^n \subset K$ and (for some $m$) there exist unit
vectors $(u_i)_1^m$ in $K$ and positive numbers $(c_i)_1^m$
satisfying:
\begin{equation*}
\sum c_iu_i=0
\end{equation*}
and
\begin{equation*}
\sum c_i\inprod{x}{u_i}^2=||x||^2 \ \text{for all} \ x \in
\mathbb{R}^n.
\end{equation*}
\end{theorem}
For the proof and discussion see \cite{ball}.
\\ \indent There also exists a unique ellipsoid of
minimal volume containing $K$, known as the Loewner ellipsoid of
$K$; we will denote it by $L_K$. It was shown by John in
\cite{john} that if $B^n$ contains $K$, then the same condition on
points in the intersection of boundaries is necessary and
sufficient for a unit ball $B^n$ to be the Loewner ellipsoid of
$K$.
\\
\begin{defin}
For a convex body $K$ in $V$ we will use $K^{\circ}$ to denote the
polar of $K$,
\begin{equation*}
K^{\circ}=\{x \in \mathbb{R}^n \mid \inprod{x}{y} \leq 1 \foral y
\in K\}.
\end{equation*}
\end{defin}
\indent The following proposition relating John and Loewner
Ellipsoids of polar bodies will be useful later on.
\begin{prop}
\label{jolo} Let $L_K$ be the Loewner ellipsoid of $K$ and suppose
that the center of $L_K$ is the origin. Then John's ellipsoid of
$K^{\circ}$ is $L_K^{\circ}$
\end{prop}
\hspace{.44cm} Now we assume that the center of the Loewner
ellipsoid of $K$ is the origin and let $\alpha$ be the coefficient
of symmetry of $K$ with respect to 0, i.e. let $\alpha$ be the
largest positive real number such that
\begin{equation*}
-\alpha K \subseteq K.
\end{equation*}
\begin{prop}
\label{symmellip} Let $K$ be a convex body and let $\alpha$ be the
coefficient of symmetry of $K$ with respect to the center of the
Loewner ellipsoid $L_K$. Then
\begin{equation*}
\sqrt{\frac{\alpha}{\text{dim} \, V}}L_K \subseteq K \subseteq
L_K.
\end{equation*}
Similarly, if $\beta$ is the coefficient of symmetry of $K$ with
respect to the center of John's ellipsoid $D_K$, then
\begin{equation*}
D_K \subseteq K \subseteq \sqrt{\frac{\beta}{\text{dim} \, V}}D_K.
\end{equation*}
\end{prop}
\begin{proof}
We will show the proposition only for the case of Loewner
ellipsoid. The other case follows by polarity. Without loss of
generality we may assume that $L_K$ is a unit ball centered at the
origin. John in \cite{john} has also shown the following:
\\For a
unit vector $v \in \mathbb{R}^n$ let $d(v)$ be the distance from
the origin of the supporting hyperplane of $K$ in the direction of
$v$:
\begin{equation*}
d(v)=\max_{x \in K} \inprod{x}{v}.
\end{equation*}
Then
\begin{equation*}
d(v)d(-v) \geq \frac{1}{\text{dim} \, V}.
\end{equation*}
Now let $w \in K$ be such that
\begin{equation*}
\inprod{v}{w}=d(v).
\end{equation*}
Since the coefficient of symmetry of $K$ is $\alpha$, it follows
that
\begin{equation*}
-\alpha w \in K \quad \text{and} \quad \inprod{-\alpha
w}{-v}=\alpha d(v).
\end{equation*}
Therefore we see that
\begin{equation*}
\alpha d(v) \leq d(-v),
\end{equation*}
and thus
\begin{equation*}
\frac{d^2(-v)}{\alpha} \geq d(v)d(-v) \geq \frac{1}{\text{dim} \,
V}.
\end{equation*}
Hence it follows that for all $v \in \mathbb{R}^n$
\begin{equation*}
d(v) \geq \sqrt{\frac{\alpha}{\text{dim} \, V}},
\end{equation*}
and therefore $K$ contains a ball of radius
$\sqrt{\frac{\alpha}{\text{dim} \, V}}$.
\end{proof}
\section{Loewner Ellipsoid of an Orbit}

Let $V$ be a finite-dimensional real vector space. Let $G$ be a
compact topological group and let $\phi : G \to GL(V)$ be a
continuous representation of $G$. There exists a $G$-invariant
probability measure $\mu$ on $G$, called the Haar measure. From
existence of Haar measure it easily follows that there exists a
$G$-invariant scalar product $\inprod{\,}{}$ that makes $\phi$
into an orthogonal representation \cite{cgroup}.
\\ Let $v \in V$ and let $O_v$ be the orbit of $v$,
\begin{equation*}
O_v=\big{\{}g(v) \mid g\in G\big{\}}.
\end{equation*}
Let $W$ denote the affine span of $O_v$,
\begin{equation*}
W=\bigg{\{}\sum \lambda_i g_i(v) \mathrel{\bigg{\arrowvert}} g_i
\in G \quad \text{and} \quad \lambda_i \in \mathbb{R} \quad
\text{such that} \quad \sum \lambda_i=1\bigg{\}},
\end{equation*}
and let $K_v$ be the convex hull of $O_v$,
\begin{equation*}
K_v=\bigg{\{}\sum \lambda_i g_i(v) \mathrel{\bigg{\arrowvert}} g_i
\in G, \lambda_i \in \mathbb{R} \quad \text{such that} \quad \sum
\lambda_i=1 \ \text{and} \ \lambda_i \geq 0 \bigg{\}}.
\end{equation*}
Since $G$ is compact, it follows that $O_v$ is compact. Therefore
$K_v$ is a full-dimensional compact convex set in $W$.
\\ \indent Let $\bar{v}$ denote the projection of $v$ into the isotypic
component of $V$ corresponding to the trivial representation:
\begin{equation*}
\bar{v}=\int_G g(v) \, d\mu.
\end{equation*}
Since $\mu$ is normalized to $1$, it follows that $\bar{v} \in W$.
Now consider the linear subspace $\overline{W}$ which is obtained
by subtracting $\bar{v}$ from $W$:
\begin{equation*}
\overline{W}=\big{\{} w-\bar{v} \mid w \in W \big{\}}.
\end{equation*}
Notice that $\overline{W}$ is a $G$-module. Since
\begin{equation*}
g(v-\bar{v})=g(v)-g(\bar{v})=g(v)-\bar{v},
\end{equation*}
it follows that $\overline{W}$ is an affine span of
$O_{v-\bar{v}}$, and $K_{v-\bar{v}}$ is $K_v-\bar{v}$. Thus
instead of $K_v$ we can consider $K_{v-\bar{v}}$, inside
$\overline{W}$.
\\ \indent Therefore we have reduced our problem to computing the Loewner
ellipsoid for a point whose orbit spans the entire space affinely.
Let $v \in V$ and, without loss of generality, assume that
\begin{equation*}
V=\text{Aff}\{O_v\}.
\end{equation*}
In this case $G$ does not fix any vector in $V$ except for the
origin. For suppose not, and let $w \in V$ be fixed by $G$. Then
\begin{equation*}
\inprod{g(v)}{w}=\inprod{g^{-1}g(v)}{g^{-1}w}=\inprod{v}{w} \
\text{for all} \ g \in G.
\end{equation*}
Therefore
\begin{equation*}
\inprod{x}{w}=\inprod{v}{w}=\text{const} \quad \text{for all}
\quad x \in \text{Aff}(O_v) \! = \! V.
\end{equation*}
Thus $w=0$.
\\Let
\begin{equation*}
V=\bigoplus_{i=1}^{k}V_i
\end{equation*}
be an orthogonal decomposition of $V$ into irreducible submodules,
and let $D_i$ be the dimension of $V_i$, with $D$ denoting the
dimension of $V$. For $x \in V$ we use $l_i(x)$ to denote
orthogonal projection of $x$ into $V_i$.
Now we prove the main theorem of this section:
\begin{theorem}
\label{orbits} The Loewner ellipsoid $L$ of $K_v$ is given by the
inequality:
\begin{equation*}
\label{lellip} \sum_{i=1}^k \frac{D_i}{\li{i}{v}}\li{i}{x} \leq D.
\end{equation*}
\end{theorem}
\begin{proof}
Since $K_v$ is the convex hull of the orbit of $v$ it follows that
$K_v$ is fixed by the action of $G$. By uniqueness of the Loewner
ellipsoid, it follows that $L$ is also fixed under the action of
$G$.
\\ \indent Now let $E$ be an ellipsoid in $V$ such that
$E$ is fixed under the action of $G$ and $v \in E$. From
invariance of $E$ under $G$ it follows that
\begin{equation*}
O_v \subseteq E,
\end{equation*}
and hence
\begin{equation*}
K_v \subseteq E.
\end{equation*}
We will minimize the volume of $E$, and then we will obtain the
Loewner ellipsoid $L$.
\\ \indent Let $w$ be the center of $E$. Since $G$
fixes $E$, it follows that $G$ also fixes $w$. But the only vector
fixed by $G$ in $V$ is the origin, and thus $w=0$. Also, from the
invariance of $E$ under the action of $G$, it follows that the
defining inequality of $E$ must have the form
\begin{equation*}
\label{gellip} \sum_{i=1}^k \lambda_i \li{i}{x} \leq 1, \quad
\text{for some} \quad \lambda_i \in \mathbb{R} \quad \text{with}
\quad \lambda_i \geq 0.
\end{equation*}
To minimize volume of $E$ we may assume that $v \in
\partial E$, or in other words
\begin{equation*}
\sum_{i=1}^k \lambda_i\li{i}{v}=1.
\end{equation*}
Also,
\begin{equation*}
\vol (E)=\vol (B^D)\prod_{i=1}^{k}\lambda_i^{-D_i}.
\end{equation*}
where $B^D$ denotes the $D$-dimensional unit ball. Thus we need to
minimize
\begin{equation*}
\label{vol} \prod_{i=1}^{k}\lambda_i^{-D_i}
\end{equation*}
subject to
\begin{equation}
\label{inbound} \sum_{i=1}^k \lambda_i\li{i}{v}=1.
\end{equation}
We apply the method of Lagrange multipliers and it follows that
\begin{equation}
\label{lagv1} cD_i\frac{\vol (E)}{\lambda_i}=\li{i}{v}, \quad
\text{for some $c \in \mathbb{R}$, and for all} \quad 1 \leq i
\leq k.
\end{equation}
Therefore
\begin{equation*}
\lambda_i\li{i}{v}=cD_i \vol (E).
\end{equation*}
We substitute this into (\ref{inbound}) and it follows that
\begin{equation*}
c=\frac{1}{\vol (E)\sum_{i=1}^{k}D_i}=\frac{1}{\vol (E)D}.
\end{equation*}
This we substitute into (\ref{lagv1}) and we see that
\begin{equation*}
\lambda_i=\frac{D_i}{\li{i}{v} D}.
\end{equation*}
Now the theorem follows.
\end{proof}
\section{Duality}
\hspace{.44cm} In this section we explicitly compute the dual cone
of the cone $C$ of nonnegative polynomials and describe some of
its properties.
\\
\begin{defin}
For $f \in P_{n,d}$ let $l_{d-2i}(f)$ denote the projection of $f$
into $H^*_{n,d-2i}$.
\end{defin}
\begin{theorem}
\label{dual} For $v \in \sph$ let $p_v \in P_{n,d}$ be as follows,
\begin{equation*}
p_v=\sum_{i=0}^{\floor{d/2}} N(n,d-2i)r^{2i}L_{n,d-2i}^v.
\end{equation*}
Then for all $f \in P_{n,d}$,
\begin{equation*}
\inprod{p_v}{f}=f(v).
\end{equation*}
\end{theorem}
\begin{proof}
We observe that
\begin{equation*}
\inprod{p_{v}}{f}=\langle \sum_{i=0}^{k} l_{d-2i}(p_v) \, , \, \,
 \sum_{i=0}^{k} l_{d-2i}(f) \rangle=\sum_{i=0}^{k}
\langle N(n,d-2i)r^{2i}L^v_{n,d-2i} \, , \,\, l_{d-2i}(f) \rangle.
\end{equation*}
Therefore it would suffice to show that for all $f \in
H^*_{n,d-2i}$,
\begin{equation*}
\langle N(n,d-2i)r^{2i}L^v_{n,d-2i} \, , \, \,f \rangle=f(v).
\end{equation*}
Let $T_v$ denote the hyperplane of all polynomials in
$H^*_{n,d-2i}$ with zero at $v$,
\begin{equation*}
T_v=\{f \in H^*_{n,d-2i} \, \mid  \, f(v)=0\}.
\end{equation*}
Since $T_v$ is a hyperplane, its orthogonal complement in
$H^*_{n,d-2i}$ is a line. Let $g \in T_v^{\perp}$. We observe that
$T_v$ is fixed by the action of $J(n,v)$. Therefore $T_v^{\perp}$
is also fixed by $J(n,v)$, and from Theorem \ref{superbasic} it
follows that $g=cr^{2i}L^v_{n,d-2i}$ for some constant $c \in
\mathbb{R}$. Let $g=N(n,d-2i)r^{2i}L^v_{n,d-2i}$. Since $g \in
T_v^{\perp}$ it follows that for all $f \in H^*_{n,d-2i}$,
\begin{equation*}
\inprod{f}{g}=cf(v) \quad \text{for some constant} \quad c \in
\mathbb{R}.
\end{equation*}
To compute $c$, we use $f=r^{2i}L^v_{n,d-2i}$ and observe that
\begin{equation*}
\inprod{r^{2i}L^v_{n,d-2i}}{g}=N(n,d-2i)\inprod{L^v_{n,d-2i}}{L^v_{n,d-2i}}.
\end{equation*}
Now from Theorem \ref{facts}, we know that
\begin{equation*}
\inprod{L^v_{n,d-2i}}{L^v_{n,d-2i}}=\frac{1}{N(n,d-2i)} \quad
\text{and} \quad L^v_{n,d-2i}=1.
\end{equation*}
Thus it follows that $c=1$ as desired.
\end{proof}
\begin{remark}
\label{evend}For even $d=2k$ we may rewrite $p_v$ as
\begin{equation*}
p_v=\sum_{i=0}^{2k} N(n,2i)r^{2k-2i}L_{n,2i}^v.
\end{equation*}
\end{remark}
\begin{cor}
\label{maxi2} Let $f \in P_{n,d}$ be such that
\begin{equation*}
\frac{||f||_{\infty}}{||f||} \geq \frac{||g||_{\infty}}{||g||}
\foral  non-zero g \in P_{n,d}.
\end{equation*}
Then $f$ is a scalar multiple of $p_{e_n}$, up to a rotation of
$\mathbb{R}^n$, and
\begin{equation*}
\frac{||f||_{\infty}}{||f||}=\sqrt{D(n,d)}.
\end{equation*}
\end{cor}
\begin{proof}
By applying an appropriate rotation and rescaling we may assume
that
\begin{equation*}
||f||_{\infty}=f(e_n)=1.
\end{equation*}
We observe that $f$ lies in the affine hyperplane $T$ of all
polynomials of integral 1 on $\sph$ and furthermore $f$ is the
shortest form on this hyperplane  by the assumption that
\begin{equation*}
\frac{||f||_{\infty}}{||f||} \geq \frac{||g||_{\infty}}{||g||}
\foral  g \in P_{n,d}.
\end{equation*}
Thus $f$ is perpendicular to $T$ and from Theorem \ref{dual} it
follows that $f$ is a multiple of $p_{e_n}$.
\end{proof}
\indent Let $C^{*}$ denote the dual cone of $C$,
\begin{equation*}
C^*=\{ f \in P_{n,2k} \mid \inprod{f}{g} \geq 0 \quad \text{for
all} \quad g \in C\}.
\end{equation*}
\begin{cor}
\label{dualcor}
$C^*$ is the conical hull of the orbit of
$p_{e_n}$, where
\begin{equation*}
p_{e_n}=\sum_{i=0}^{2k} N(n,2i)r^{2k-2i}L_{n,2i}.
\end{equation*}
\end{cor}
\begin{proof}
Let $K$ be the conical hull of the points $p_v$ for all $v \in
\sph$,
\begin{equation*}
K=\biggl{\{}\sum_i \lambda_i p_{v_i} \mathrel{\bigg{\arrowvert}}
v_i \in \sph \quad \text{and} \quad \lambda_i \in \mathbb{R},
\quad \lambda_i \geq 0 \biggr{\}}.
\end{equation*}
Consider $K^*$,
\begin{equation*}
K^*=\{f \in P_{n,2k} \mid \inprod{f}{g} \geq 0 \quad \text{for
all} \quad g \in K \}
\end{equation*}
\begin{equation*}
=  \{f \in P_{n,2k} \mid \inprod{f}{p_v} \geq 0 \quad \text{for
all} \quad v \in \sph \}.
\end{equation*}
From Theorem \ref{dual}, we know that $\inprod{f}{p_v}=f(v)$, and
therefore $K^*=C$. Since $C$ is a closed cone, by the BiPolar
Theorem it follows that $K=C^*$.
\\ \indent Now let $A \in SO(n)$ be such that $Aw=v$. Then we note that
$AL^v_{n,2i}=L^w_{n,2i}$ and therefore
\begin{equation*}
Ap_v=p_w.
\end{equation*}
Thus the set of $p_v$ for all $v \in \sph$ is the same as the
orbit of $p_{e_n}$ and we obtain the desired result.
\end{proof}

\section{John's Ellipsoid of the Cone of Nonnegative Polynomials}
\hspace{.44cm} In this section we compute John's Ellipsoid of
$\overline{C}_{n,2k}$. Recall that $M$ is the hyperplane of all
forms of integral 1 on $\sph$. If we regard the point
$r^{2k}=(x_1^2+\ldots+x_n^2)^k$ as the origin in $M$, then the
inner product $\inprod{}{}$ induces an inner product in $M$ which
we will denote $\inprodm{}{}$,
\begin{equation*}
\inprodm{f}{g}=\inprod{f-r^{2k}}{g-r^{2k}} \quad \text{for} \quad
f,g \in M.
\end{equation*}
Recall that $C^*$ is the dual cone of $C$, and define
$\overline{C}^*$ by intersecting $C^*$ with the hyperplane $M$,
\begin{equation*}
\overline{C}^*=C^* \cap M.
\end{equation*}
We now establish a relationship between $\overline{C}$ and
$\overline{C}^*$ in terms of $\inprodm{}{}$.

\begin{lemma}
\label{convert} Let $\overline{C}^{\circ}$ be the polar of
$\overline{C}$ with respect to $\inprodm{}{}$. Then
\begin{equation*}
\overline{C}^{\circ}=-\overline{C}^*+2r^{2k}.
\end{equation*}
\end{lemma}
\begin{proof}
We observe that
\begin{equation*}
\overline{C}^{\circ}=\{f \in M \mid \inprodm{f}{g} \leq 1 \foral g
\in \overline{C}\}
\end{equation*}
\begin{equation*}
=\{f \in M \mid \inprod{f-r^{2k}}{g-r^{2k}} \leq 1 \foral g \in
\overline{C}\}.
\end{equation*}
Since both $f$ and $g$ have integral 1 on $\sph$, it follows that
\begin{equation*}
\inprod{f}{r^{2k}}=\inprod{g}{r^{2k}}=1,
\end{equation*}
and therefore
\begin{equation*}
\overline{C}^{\circ}=\{f \in M \mid \inprod{f}{g} \leq 2 \foral g
\in \overline{C}\}.
\end{equation*}
Thus
\begin{equation*}
-\overline{C}^{\circ}+2r^{2k}=\{f \in M \mid \inprod{f}{g} \geq 0
\foral g \in \overline{C}\}=\overline{C}^*.
\end{equation*}
\end{proof}
\begin{theorem}
\label{firstellipsoid} The Loewner Ellipsoid $E$ of
$\overline{C}^*$ is a ball with center $r^{2k}$ and radius
\begin{equation*}
\sqrt{D(n,2k)-1}=\sqrt{\binom{n+2k-1}{2k}-1}
\end{equation*}
\end{theorem}
\begin{proof}
From Corollary \ref{dualcor} it follows that $\overline{C}^*$ is
the convex hull of the orbit of $p_{e_n}$. Therefore we can apply
Theorem \ref{orbits}. The irreducible subspaces are $H_{n,2i}^*$
for $1 \leq i \leq k$. Let $l_{2i}(f)$ denote the projection of
$f$ into $H^*_{n,2i}$ and then
\begin{equation*}
||l_{2i}(p_{e_n})||^2=||N(n,2i)L_{n,2i}||^2=N(n,2i).
\end{equation*}
The result now follows from Theorem \ref{orbits}.
\end{proof}

\begin{theorem}
\label{main} John's ellipsoid $D$ of $\overline{C}$ is a ball with
center $r^{2k}$ and radius
\begin{equation*}
\frac{1}{\sqrt{D(n,2k)-1}}=\frac{1}{\sqrt{\binom{n+2k-1}{2k}-1}}
\end{equation*}
\end{theorem}
\begin{proof}
From Lemma \ref{convert} we know that
\begin{equation*}
\overline{C}^{\circ}=-\overline{C}^*+2r^{2k}.
\end{equation*}
Therefore the Loewner ellipsoid of $\overline{C}^{\circ}$ is a
ball with center $r^{2k}$ and radius $\sqrt{D(n,2k)}$. By
Proposition \ref{jolo} we know that John's Ellipsoid of
$\overline{C}$ is the polar of the Loewner ellipsoid of
$\overline{C}^{\circ}$ and the theorem follows.
\end{proof}

\section{Coefficient of Symmetry of The Cone of Nonnegative Polynomials}
\hspace{.44cm} In this section we compute the coefficient of
symmetry of $\overline{C}$ with respect to $r^{2k}$. We begin by
showing that all forms symmetric with respect to an axis are sums
of squares of forms.\\
\begin{defin}
For $0 \leq a \leq 1$, let
\begin{equation*}
q_a(x)=x_n^2-ar^{2}=(1-a)x_n^2-a(x_1^2+\ldots +x_{n-1}^2).
\end{equation*}
\end{defin}
\begin{defin}
For $f \in P_{n,d}$ let $V(f)$ be the vanishing set of $f$,
\begin{equation*}
V(f)=\{x \in \mathbb{R}^n \mid f(x)=0\}.
\end{equation*}
\end{defin}
\begin{lemma}
\label{squares} Let $f(x) \in P_{n,2k}$ be a nonnegative form and
suppose that $f$ is fixed by $J(n,v)$ for some $v\in
\mathbb{R}^n$. Then $f$ is a sum of squares of forms.
\end{lemma}
\begin{proof}
We induct on $k$.
\\\emph{Base Case}: $k=1.$ In this case we are dealing with
homogeneous quadratics and all nonnegative homogeneous quadratics
are sums of squares.
\\\emph{Inductive step}: $k \Rightarrow k+1.$ Applying a suitable
rotation of $\mathbb{R}^n$, we may assume that $f$ is fixed by
$J(n,e_n)$. It will suffice to show the lemma for $f$ with a zero,
since we can consider the form
\begin{equation*}
f-\alpha r^{2k},
\end{equation*}
where $\alpha$ is the minimum of $f$ on $\sph$. Since $f$ has a
zero, and $f$ is fixed by $J(n,e_n)$ it follows $V(f)$ is a
nonempty subset of $\mathbb{R}^n$ that is fixed by $J(n,e_n)$.
Hence $V(f)$ contains $V_{q_a}$ for some $a \in [0,1]$.
\\ We first deal with the two degenerate cases:
\\ \indent If $a=1$ then
\begin{equation*}
q_1=-(x_1^2+\ldots +x_{n-1}^2) \quad \text{and} \quad V_{q_a}=\{
\lambda e_n \mid \lambda \in \mathbb{R} \}
\end{equation*}
Since $f(e_n)=0$ we can write
\begin{equation*}
f=\sum_{i=0}^{2k-1} x_n^i g_i.
\end{equation*}
where $g_i$ depend only on $x_1, \ldots, x_{n-1}$. Since $f$ is
fixed by $J(n,e_n)$, it follows that $g_{i}$ is fixed by
$J(n,e_n)$ for all $0 \leq i  \leq 2k-1$. Since $g_{i}$ depends
only on $x_1, \ldots, x_{n-1}$, we see that $g_{i}$ is fixed by
$SO(n-1)$. Then $i$ must be even and
\begin{equation*}
g_{2i}=\lambda_i(x_1^2+ \cdots + x_{n-1}^2)^{k-i} \quad \text{for
some} \quad \lambda_i \in \mathbb{R}.
\end{equation*}
Thus $x_1^2+ \cdots + x_{n-1}^2$ divides $f$. We write
$f=(x_1^2+\cdots +x_{n-1}^2)g$ and $g$ is sum of squares by
induction, and then so is $f$.
\\ \indent If $a=0$, then $q_0(x)=x_n^2$ and $x_n$ divides $f$,
but since $f$ is nonnegative, it follows that $x_n^2$ divides $f$
and $f=x_n^2g$. By induction, $g$ is a sum of squares, and then
$f$ is as well .
\\ \indent For $0< a < 1 $, let $I=I(V_{q_a})$ be the vanishing ideal of
$V_{q_a}$:
\begin{equation*}
I=\{f \in \mathbb{R}[x_1, \ldots ,x_n] \mid f(x)=0 \quad \text{for
all} \quad x \in V_{q_a}\},
\end{equation*}
where $\mathbb{R}[x_1, \ldots ,x_n]$ is the ring of real
polynomials in $n$ variables. We will show that $I$ is a principal
ideal generated by $q_a$.
\\\indent Let $g \in I$. By reducing modulo $q_a$ we may write
\begin{equation*}
g=bq_a + x_nc+d,
\end{equation*}
where $c$ and $d$ are polynomials that depend only on $x_1,
\ldots, x_{n-1}$. Let
\begin{equation*}
h=g-bp_a=x_nc+d.
\end{equation*}
We observe that $h \in I$ and also $h(x_1,\ldots,x_{n-1},-x_n) \in
I$, since $V(q_a)$ is fixed by reflection about the $e_1, \ldots
,e_{n-1}$ hyperplane. Thus $-x_nc+d \in I$, and then $x_nc$ and
$d$ are in $I$. But since $a>0$, the vanishing set of $q_a$
intersects the hyperplane $x_n=0$ only at the origin. Thus we see
that $c \in I$. Also, $c$ and $d$ only depend on the first $n-1$
variables. Therefore, since $a<1$, we see that $c$ and $d$ vanish
on the entire hyperplane spanned by $e_1, \ldots, e_{n-1}$. Hence,
\begin{equation*}
c=d \equiv 0.
\end{equation*}
Thus $I=(q_a)$.
\\ \indent Since $I$ is a principal ideal generated by $q_a$ and $f \in
I$ it follows that $q_a$ divides $f$, and we can write $f=q_ag$.
Now we note that $q_a(x) \geq 0 \ \text{for} \ x \in \sph \,
\text{with} \ x_n^2 \geq a, \ \text{and} \ q_a(x) < 0 \ \text{for}
\ x \in \sph \ \text{with} \ x_n^2 < a$. Since $q_ag \geq 0$, it
follows that
\begin{equation*}
g(x)=0 \ \text{for all} \ x \in V(p_a),
\end{equation*}
otherwise the sign of $g$ does not change in the neighborhood of
some $x \in V(q_a)$, which yields a contradiction since $a < 1$.
Thus $g \in I$ and therefore $q_a$ divides $g$. Hence $q_a^2$
divides $f$. We write $f=q_a^2h$ and $h$ is a sum of squares by
induction.
\end{proof}
\begin{remark}
\label{factor} From the proof of Lemma \ref{squares} it follows
that if
\begin{equation*}
V(q_a) \subseteq V(f) \quad \text{with} \quad 0 \leq a <1
\end{equation*}
for some nonnegative $f \in P_{n,2k}$, not necessarily symmetric
with respect to $J(n,e_n)$, then
\begin{equation*}
q_a^2 \quad \text{divides} \quad f \quad \text{for} \quad 0<a<1
\quad \text{and} \quad x_n^2 \quad \text{divides} \quad f \quad
\text{if} \quad a=0.
\end{equation*}
\end{remark}
Our goal is to compute the coefficient of symmetry of
$\overline{C}$. We begin with the crucial integral inequality. \\
\begin{defin} Let $\text{Max}$ denote the maximal $L^{\infty}$ norm for the
functions in $\overline{C}$,
\begin{equation*}
\text{Max}=\max_{f \in \overline{C}} ||f||_{\infty}.
\end{equation*}
\end{defin}
\begin{theorem}
\label{maxi3} Let $f \in \overline{C}$ be such that
$||f||_{\infty}=\text{Max}$. Then
\begin{equation*}
f=\frac{1}{D(n,k)}\left(\sum_{l=0}^{\lfloor k/2 \rfloor
}N(n,k-2l)r^{2l}L_{n,k-2l}\right)^2,
\end{equation*}
up to a rotation of $\mathbb{R}^n$, and
\begin{equation*}
\text{Max}=||f||_{\infty}=D(n,k).
\end{equation*}
\end{theorem}
\begin{proof}
Let $f \in \overline{C}$ be such that $||f||_{\infty}=\text{Max}$.
Applying a rotation of $\mathbb{R}^n$, if necessary, we may assume
that $f(e_n)=\text{Max}$. Now let $p$ be the average of $f$ over
$J(n,e_n)$,
\begin{equation*}
p=\int_{A \in J(n,e_n)}Af \, d\mu,
\end{equation*}
where $\mu$ is the normalized Haar measure on $J(n,e_n)$. Clearly,
$p$ is a nonnegative form and
\begin{equation*}
\int_{S^{n-1}}p \, d\sigma=1.
\end{equation*}
Thus $p \in \overline{C}$. Also,
$||p||_{\infty}=||f||_{\infty}=\text{Max}$, since $p(e_n)=f(e_n)$.
Since $p$ is the average of $f$ over $J(n,e_n)$, it follows that
$p$ is fixed by $J(n,e_n)$. Then from Lemma \ref{squares} we see
that $p$ is a sum of squares.
\\ \indent Since $p \in \overline{Sq}$, it is a convex combination of extreme
points of $\overline{Sq}$, and an extreme point of $\overline{Sq}$
must be a square. Thus we see that
\begin{equation*}
p=\sum \lambda_i h_i^2 \quad \text{with} \quad \lambda_i > 0,
\quad \sum \lambda_i=1,
\end{equation*}
where $h_i \in P_{n,k}$
Therefore,
\begin{equation}
\label{pepsi} \text{Max}=p(e_n)=\sum \lambda_i h_i^2(e_n).
\end{equation}
But
\begin{equation*}
\text{Max} \geq ||h_i^2||_{\infty} \geq h_i^2(e_n) \quad \text{and
therefore} \quad ||h_i^2||_{\infty}=h_i^2(e_n)=\text{Max}.
\end{equation*}
Thus there exists $h \in P_{n,k}$ such that $h^2 \in
\overline{Sq}$ and
\begin{equation*}
||h^2||_{\infty}=\text{Max}.
\end{equation*}
Then we observe that
\begin{equation*}
\frac{||h||_{\infty}}{||h||} \geq \frac{||g||_{\infty}}{||g||}
\foral g \in P_{n,k}.
\end{equation*}
Then from Corollary \ref{maxi2} it follows that
\begin{equation*}
h=\frac{1}{\sqrt{D(n,k)}}\sum_{i=0}^{\floor{k/2}}N(n,k-2i)r^{2i}L_{n,k-2i},
\end{equation*}
up to a rotation of $\mathbb{R}^n$. Also from Corollary
\ref{maxi2} we know that
\begin{equation*}
\text{Max}=||h^2||_{\infty}=D(n,k).
\end{equation*}
\indent Now we will show that up to a rotation of $\mathbb{R}^n$,
the only form in  $\overline{C}$ with maximal $L^{\infty}$ norm is
$h^2$. We know that all Legendre harmonics are fixed by
$J(n,e_n)$. Therefore it follows that $h^2$ is also fixed by
$J(n,e_n)$. Now we observe that from the proof of Lemma
\ref{maxi2} it is clear that
\begin{equation*}
h^2(e_n)=h^2(-e_n)=||h^2||_{\infty}=\text{Max},
\end{equation*}
and $e_n$, $-e_n$ are the only points where the maximum occurs.
Thus, if $A \in SO(n)$ acts on $h$, then it ether fixes $h$, or
the maximum of $Ah$ occurs not at $\pm e_n$. Therefore $h^2$ is
the only square, and thus the only extreme point of
$\overline{Sq}$, which takes on the value $\text{Max}$ at $e_n$.
Now going back to \eqref{pepsi} we see that $p=h^2$, since $p$ is
a convex linear combination of extreme points of $\overline{Sq}$
with value $\text{Max}$ at $e_n$. Therefore
\begin{equation*}
h^2=\int_{A\in J(n,\epsilon_n)}Af \, d\mu.
\end{equation*}
Now $h^2$ is lies in the boundary of $\overline{C}$, and thus it
must have a zero. But $h^2$ is also the average of $f$ over
$J(n,e_n)$ and we know that $f$ is nonnegative. Therefore we see
that $V(f)$ contains $V(q^2)$. Since $q^2$ is fixed by $J(n,e_n)$,
it follows that $V(q_a) \subseteq V(h^2)$, for some $a \in [0,1]$,
and since $h(e_n) \neq 0$ it follows that $0 \leq a <1$. Then it
follows from Remark \ref{factor} that we can factor out a square
of a form fixed by $J(n,e_n)$ from $h^2$ and $f$. Call it $m^2$,
and let
\begin{equation*}
\tilde{h}^2=\frac{h^2}{m^2} \quad \text{and} \quad
\tilde{f}=\frac{f}{m^2}.
\end{equation*}
Again, $\tilde{h}^2$ generates an extreme ray of a cone of sums of
squares, now of a lesser degree, otherwise $h^2$ would not lie on
an extreme ray. Since $m$ is fixed by $J(n,e_n)$, we still have
\begin{equation*}
\tilde{h}^2=\int_{A\in J(n,\epsilon_n)}A\tilde{f}d\mu \, ,
\end{equation*}
because averaging over $J(n,e_n)$ is the same as taking the
average over slices of the sphere with hyperplanes perpendicular
to $e_n$. Thus, again by proof of Lemma \ref{squares} we can
factor out the same square from both $\tilde{h}^2$ and $\tilde{f}$
and we can continue with this process, and in the end $f=h^2$.
\end{proof}
\begin{cor}
\label{ineqcor1} For all $f \in P_{n,k}$
\begin{equation*}
||f||_{\infty} \leq {\binom{n+kl-1}{kl}}^{\frac{1}{2l}}||f||_{2l}.
\end{equation*}
\end{cor}
\begin{proof}
We apply Theorem \ref{maxi3} to $f^{2l}$. Since $f^{2l}$ is
nonnegative, and $f^{2l} \in P_{n,2kl}$, from Theorem \ref{maxi3}
we know that
\begin{equation*}
||f^{2l}||_{\infty} \leq D(n,kl) ||f^{2l}||_1.
\end{equation*}
Since
\begin{equation*}
||f^{2l}||_{\infty}=||f||_{\infty}^{2l} \quad \, \text{and} \quad
\, ||f^{2l}||_1= ||f||_{2l}^{2l},
\end{equation*}
by taking $2l$-th root of both sides we obtain the desired
inequality.
\end{proof}
\begin{theorem}
\label{nonnegcoeffsym} The coefficient of symmetry of
$\overline{C}$ with respect to $r^{2k}$ is
\begin{equation*}
\frac{1}{D(n,k)-1}=\frac{1}{\binom{n+k-1}{k}-1}.
\end{equation*}
\end{theorem}
\begin{proof}
Let $f \in \partial \overline{C}$, and denote by $\bar{f}$ the
polynomial in $\partial \overline{C}$ that is opposite to $f$ with
respect to $r^{2k}$,
\begin{equation*}
\label{fbar} \bar{f}=\alpha(r^{2k}-f)+r^{2k} \quad \text{for some}
\quad \alpha \in \mathbb{R} \quad \text{such that} \quad \alpha >
0.
\end{equation*}
Since $\bar{f} \in \partial \overline{C}$, it is a nonnegative
form with a zero. Then it follows that
\begin{equation*}
\label{coeffsym} \alpha =\frac {1}{\max_{x \in \sph} f(x)-1}=
\frac {1}{||f||_{\infty}-1}.
\end{equation*}
Thus
\begin{equation}
\label{fbarfin} \bar{f}= \frac
{1}{||f||_{\infty}-1}(r^{2k}-f)+r^{2k},
\end{equation}
and, since the minimum of $f$ on $\sph$ is zero,
\begin{equation*}
||\bar{f}||_{\infty}=\alpha+1=\frac{||f||_{\infty}}{||f||_{\infty}-1}.
\end{equation*}
Also using (\ref{fbarfin}) we see that:
\begin{equation*}
\frac{||f-r^{2k}||}{||\bar{f}-r^{2k}||}=\frac{||f-r^{2k}||}{||\frac
{1}{||f||_{\infty}-1}(r^{2k}-f)||}=||f||_{\infty}-1.
\end{equation*}
Therefore it follows that the coefficient of symmetry of
$\overline{C}$ with respect to $r^{2k}$ is
$\frac{1}{\text{Max}-1}$. From Theorem \ref{maxi3}, we know that
$\text{Max}=D(n,k)$, and the result follows.
\end{proof}
\begin{cor}
\label{ineqcor} Let $M_f$ denote the maximum of $f$ on $\sph$ and
let $m_f$ denote the minimum of $f$ on $\sph$. Let
$\alpha=\frac{1}{\binom{n+k-1}{k}}$. Then,
\begin{equation*}
\alpha M_f +(1-\alpha)m_f \leq \int_{\sph} f \, d\sigma \leq
(1-\alpha) M_f +\alpha m_f,
\end{equation*}
and both inequalities are sharp.
\end{cor}
\begin{proof}
Consider the set $W$ obtained from $\overline{C}$ by subtracting
$r^{2k}$ from all forms in $\overline{C}$,
\begin{equation*}
W=\overline{C}-r^{2k}.
\end{equation*}
We observe that $W$ is the set of all forms of integral zero with
minimum at most $-1$ on $\sph$.
\\\indent From the definition of
$W$ it follows that the coefficient of symmetry of $W$ around 0 is
the same as the coefficient of symmetry of $\overline{C}$ around
$r^{2k}$. Thus the coefficient of symmetry of $W$ around 0 is,
\begin{equation*}
\frac{1}{\binom{n+k-1}{k}-1}=\frac{\alpha}{1-\alpha}.
\end{equation*}
But since
\begin{equation*}
M_{-f}=-m_{f} \quad \text{and} \quad m_{-f}=-M_f,
\end{equation*}
it follows that for all $f \in P_{n,2k}$ of integral 0,
\begin{equation*}
\frac{\alpha}{1-\alpha} \leq \frac{-M_f}{m_f} \leq
\frac{1-\alpha}{\alpha}.
\end{equation*}
For $f \in P_{n,2k}$ consider
\begin{equation*}
\hat{f}=f-\left( \int_{\sph} f \, d \sigma \right) r^{2k}.
\end{equation*}
We have shown above that
\begin{equation*}
\frac{-(M_f-\int_{\sph} f \, d \sigma)}{m_f-\int_{\sph} f \, d
\sigma}=\frac{-M_{\hat{f}}}{m_{\hat{f}}} \leq
\frac{1-\alpha}{\alpha}.
\end{equation*}
Thus
\begin{equation*}
\frac{1}{\alpha} \int_{\sph} f \, d \sigma \ \geq \ M_f +
\left(\frac{1}{\alpha}-1\right)m_f,
\end{equation*}
and one side of the desired inequality follows. The other half is
done in the same way.
\end{proof}
\begin{cor}
\label{mainint} $\overline{C}$ is contained in ball of radius
\begin{equation*}
\sqrt{D(n,k)-1}=\sqrt{\binom{n+k-1}{k}-1} \, ,
\end{equation*}
or, equivalently, for all nonnegative $f \in P_{n,2k}$
\begin{equation*}
||f||_2 \leq \sqrt{\binom{n+k-1}{k}}||f||_1.
\end{equation*}
\end{cor}
\begin{proof}
From Theorem \ref{main} we know that John's ellipsoid of
$\overline{C}$ is a ball of radius
$\frac{1}{\sqrt{\smash[b]{D(n,2k)-1}}}$ around $r^{2k}$, and the
coefficient of symmetry of $\overline{C}$ with respect to $r^{2k}$
is $\frac{1}{D(n,k)-1}$. We apply Proposition \ref{symmellip}, and
it follows that therefore $\overline{C}$ is contained in the ball
of radius
\begin{equation*}
\frac{\sqrt{D(n,k)-1}\sqrt{D(n,2k)-1}}{\sqrt{D(n,2k)-1}}={\sqrt{D(n,k)-1}}
=\sqrt{\binom{n+k-1}{k}-1},
\end{equation*}
centered at $r^{2k}$, as desired.
\end{proof}
\section{Cone of Sums of Powers of Linear Forms}
\hspace{.44cm} In order to study the cone $\lf$ we will need to
decompose $x_n^{2k}$ as a sum of Legendre harmonics. We begin by
recalling the Rodrigues rule.
\begin{lemma}
Rodrigues Rule, \cite{mull}: Let $Q_{n,d}(t)$ be the Legendre
polynomial defined in the Preliminaries. Then
\begin{equation*}
\label{rodr}
\int^{+1}_{-1}f(t)Q_{n,d}(t)(1-t^2)^{\frac{n-3}{2}}dt=
R_d(n)\int_{-1}^{+1}f^{(n)}(t)(1-t^2)^{\frac{2d+n-3}{2}},
\end{equation*}
where $R_d(n)$ is the Rodrigues constant:
\begin{equation*}
\label{rodc}
R_d(n)=\frac{\Gamma(\frac{n-1}{2})}{2^n\Gamma(\frac{4k+n-1}{2})}.
\end{equation*}
\end{lemma}
Since $x_n^{2k}$ is symmetric is fixed by the action of $J(n,e_n)$
it decomposes as a sum of the Legendre harmonics. The next theorem
gives the precise decomposition.
\begin{theorem}
\label{powerdecomp}
\begin{equation*}
\frac{x_n^{2k}}{\int_{S^{n-1}}x_n^{2k} \, d\sigma}=
\sum_{l=0}^{k}\frac{k!\Gamma(\frac{2k+n}{2})}
{(k-l)!\Gamma(\frac{2k+2l+n}{2})}N(n,2l)r^{2k-2l}L_{n,2l}.
\end{equation*}
\end{theorem}
\begin{proof}
We first recall the well-known fact that
\begin{equation}
\label{basic} \int_{S^{n-1}}x_n^{2k}d\sigma=
\frac{\Gamma(\frac{2k+1}{2})\Gamma(\frac{n}{2})}{\sqrt{\pi}\Gamma(\frac{n+2k}{2})}.
\end{equation}
See, for example, \cite{barv}.
\\Since $x_n^{2k}$ is fixed by $J(n,e_n)$, we know that it decomposes
as a sum of Legendre harmonics of even degrees. Therefore it
suffices to compute
\begin{equation*}
\inprod{x_n^{2k}}{r^{2l}L_{n,2k-2l}}
=\int_{\sph}x_n^{2k}L_{n,2k-2l} \, d\sigma.
\end{equation*}
On $S^{n-1}$ both $x_n^{2k}$ and $L_{n,2k-2l}$ are functions of
the last coordinate, and hence this integral translates into
\begin{equation*}
\frac{|S^{n-2}|}{|S^{n-1}|}
\int^{+1}_{-1}t^{2k}Q_{n,2k-2l}(t)(1-t^2)^{\frac{n-3}{2}}dt,
\end{equation*}
where $|\sph|$ denotes the surface area of $\sph$. Now we apply
the Rodrigues Rule to
\begin{equation*}
\int^{+1}_{-1}t^{2k}Q_{n,2k-2l}(t)(1-t^2)^{\frac{n-3}{2}}dt,
\end{equation*}
and get:
\begin{equation*}
\frac{(2k)!}{(2l)!}R_{2k-2l}(n)\int^{+1}_{-1}t^{2l}(1-t^2)^{\frac{4k-4l+n-3}{2}}dt.
\end{equation*}
This we can interpret back as an integral over the sphere of
dimension $4k-4l+n-1$ and we obtain:
\begin{equation*}
\frac{(2k)!|S^{4k-4l+n-1}|}{(2l)!|S^{4k-4l+n-2}|}R_{2k-2l}(n)
\int_{S^{4k-4l+n-1}}x_n^{2l}d\sigma.
\end{equation*}
Next we substitute in (\ref{basic}) to get,
\begin{equation*}
\frac{(2k)!|S^{4k-4l+n-1}|\Gamma(\frac{2l+1}{2})\Gamma(\frac{4k-4l+n}{2})}
{\sqrt{\pi}(2l)!|S^{4k-4l+n-2}|
\Gamma(\frac{4k-2l+n}{2})}R_{2k-2l}(n).
\end{equation*}
Now,
\begin{equation*}
|S^{n-1}|=\frac{\pi^{\frac{n}{2}}}{\Gamma(\frac{n}{2})},
\end{equation*}
and thus we get,
\begin{equation*}
\frac{(2k)!\Gamma(\frac{2l+1}{2})\Gamma(\frac{4k-4l+n-1}{2})}
{(2l)!\Gamma(\frac{4k-2l+n}{2})}R_{2k-2l}(n).
\end{equation*}
Substituting in the value of $R_{2k-2l}(n)$ we obtain
\begin{equation*}
\int^{+1}_{-1}t^{2k}Q_{n,2k-2l}(t)(1-t^2)^{\frac{n-3}{2}}dt=
\frac{(2k)!\Gamma(\frac{2l+1}{2})\Gamma(\frac{n-1}{2})}
{2^{2k-2l}(2l)!\Gamma(\frac{4k-2l+n}{2})}.
\end{equation*}
Thus we get that
\begin{equation*}
\inprod{x_n^{2k}}{r^{2l}L_{n,2k-2l}}=
\frac{|S^{n-2}|(2k)!\Gamma(\frac{2l+1}{2})\Gamma(\frac{n-1}{2})}
{|S^{n-1}|2^{2k-2l}(2l)!\Gamma(\frac{4k-2l+n}{2})}=
\end{equation*}
\begin{equation}
\label{semi} \frac{(2k)!\Gamma(\frac{2l+1}{2})\Gamma(\frac{n}{2})}
{\sqrt{\pi}2^{2k-2l}(2l)!\Gamma(\frac{4k-2l+n}{2})}.
\end{equation}
Now the doubling rule for Gamma function says that
\begin{equation*}
2^{x-1}\Gamma\left(\frac{x}{2}\right)\Gamma\left(\frac{x+1}{2}\right)
=\sqrt{\pi}\Gamma(x).
\end{equation*}
Applying this to $x=2l+1$ we get
\begin{equation}
\label{dbl}
\Gamma\left(\frac{2l+1}{2}\right)=\frac{\sqrt{\pi}\Gamma(2l+1)}{2^{2l}\Gamma(l+1)}.
\end{equation}
Substituting \eqref{dbl} into \eqref{semi} we have
\begin{equation*}
\label{fin} \langle x_n^{2k},r^{2l}L_{n,2k-2l} \rangle
=\frac{(2k)!\Gamma(\frac{n}{2})}
{l!2^{2k}\Gamma(\frac{4k-2l+n}{2})}.
\end{equation*}
Thus using \eqref{basic},
\begin{equation}
\label{final} \frac{\langle x_n^{2k},r^{2l}L_{n,2k-2l} \rangle
}{\int_{S^{n-1}}x_n^{2k}d\sigma}=
\frac{\sqrt{\pi}(2k)!\Gamma(\frac{2k+n}{2})}
{2^{2k}l!\Gamma(\frac{4k-2l+n}{2})\Gamma(\frac{2k+1}{2})}.
\end{equation}
Now we again apply the doubling rule this time to $x=2k+1$ to get
\begin{equation*}
\Gamma
\left(\frac{2k+1}{2}\right)=\frac{\sqrt{\pi}\Gamma(2k+1)}{2^{2k}\Gamma(k+1)},
\end{equation*}
which we substitute into \eqref{final}:
\begin{equation*}
\label{finally}
\frac{\inprod{x_n^{2k}}{r^{2l}L_{n,2k-2l}}}{\int_{\sph} x_n^{2k}
\, d\sigma}= \frac{k!\Gamma(\frac{2k+n}{2})}
{l!\Gamma(\frac{4k-2l+n}{2})}.
\end{equation*}
Now recall that
\begin{equation*}
||r^{2l}L_{n,2k-2l}||^2=||L_{n,2k-2l}||^2=\frac{1}{N(n,2k-2l)},
\end{equation*}
and the desired result follows.
\end{proof}
\indent We now make a crucial definition.
\\
\begin{defin}
Let $T_{2m,2k}:P_{n,2k} \to P_{n,2k}$ be a linear operator defined
by
\begin{equation*}
(T_{2m,2k}p)(x)=\frac{\int_{\sph} p(v) \inprod{x}{v}^{2m} \,
d\sigma(v)}{\int_{\sph}x_n^{2m} \, d\sigma} \quad \text{for} \quad
x,v \in \sph.
\end{equation*}
\end{defin}
\begin{remark}
It will follow from Lemma \ref{rewrite} that the operators
$T_{2m,2k}$ have been defined in a different form by Reznick in
$\cite{rez1}$.
\end{remark}
We observe that $T_{2m,2k}$ maps nonnegative forms to the sums of
powers of linear forms. The following lemma shows the precise
action of $T_{2m,2k}$ on $P_{n,2k}$.
\begin{lemma}
\label{rewrite}
\begin{equation*}
T_{2m,2k}(f)=\sum_{i=0}^{k}\frac{m!\Gamma(\frac{2m+n}{2})}
{(m-i)!\Gamma(\frac{2m+2i+n}{2})}l_{2i}(f).
\end{equation*}
\end{lemma}
\begin{proof}
We rewrite $T_{2m,2k}f$ as
\begin{equation*}
(T_{2m,2k}f)(x)=\inprod{f}{\frac{\inprod{x}{v}^{2m}}{\int_{\sph}\inprod{x}{v}^{2m}
\, d\sigma(v)}}
\end{equation*}
We apply Theorem \ref{powerdecomp} and it follows that
\begin{equation*}
(T_{2m,2k}f)(x)=\inprod{f}{\sum_{i=0}^{m}\frac{m!\Gamma(\frac{2m+n}{2})}
{(m-i)!\Gamma(\frac{2m+2i+n}{2})}N(n,2i)L^v_{n,2i}}.
\end{equation*}
Now we decompose $f$ as a sum of spherical harmonics and observe
that
\begin{equation*}
(T_{2m,2k}f)(x)=\sum_{i=0}^{k}
\inprod{l_{2i}(f)}{\frac{m!\Gamma(\frac{2m+n}{2})}
{(m-i)!\Gamma(\frac{2m+2i+n}{2})}N(n,2i)L^v_{n,2i}}.
\end{equation*}
We recall that by Theorem \ref{dual}
\begin{equation*}
\inprod{l_{2i}(f)}{N_{n,2i}L_{n,2i}^v}=(l_{2i}f)(v).
\end{equation*}
Therefore
\begin{equation*}
T_{2m,2k}(f)=\sum_{i=0}^{k}\frac{m!\Gamma(\frac{2m+n}{2})}
{(m-i)!\Gamma(\frac{2m+2i+n}{2})}l_{2i}(f).
\end{equation*}
\end{proof}
\begin{remark}
\label{diagonal} It follows from Lemma \ref{rewrite} that
$T_{2m,2k}$ is a diagonal operator on the harmonic subspaces of
$P_{n,2k}$. Thus $T_{2m,2k}$ commutes with the action of $SO(n)$.
\end{remark}
\begin{theorem}
\label{mainpowers} The Loewner Ellipsoid of $\overline{\lf}$ is
given by the inequality
\begin{equation*}
\sum_{i=1}^{k} \left(\frac{(k-i)!\Gamma(\frac{2k+2i+n}{2})}
{k!\Gamma(\frac{2k+n}{2})}\right)^2 \li{2i}{f} \leq D_{n,2k}-1,
\end{equation*}
and the coefficient of symmetry of $\overline{\lf}$ is
\begin{equation*}
\frac{1}{D(n,k)-1}=\frac{1}{\binom{n+k-1}{k}-1}.
\end{equation*}
\end{theorem}
\begin{proof}
By Lemma \ref{rewrite} and Theorem \ref{powerdecomp},
\begin{equation*}
T_{2k,2k}(p_{e_n})=\sum_{i=0}^{k}\frac{k!\Gamma(\frac{2k+n}{2})}
{(k-i)!\Gamma(\frac{2k+2i+n}{2})}N(n,2i)L_{n,2i}=\frac{x_n^{2k}}{\int_{\sph}
x_n^{2k} \, d\sigma}.
\end{equation*}
Since $T_{2k,2k}$ commutes with the action of $SO(n)$, it follows
that
\begin{equation*}
T_{2k,2k}(\overline{C}^*)=\overline{\lf}.
\end{equation*}
Therefore $T_{2k,2k}$ maps the Loewner ellipsoid of
$\overline{C}^*$ to the Loewner ellipsoid of $\overline{\lf}$. By
Theorem \ref{firstellipsoid} the Loewner ellipsoid of
$\overline{\lf}$ is a ball with center $r^{2k}$ and of radius
$\sqrt{D(n,k)-1}$. The inequality for the Loewner ellipsoid of
$\overline{\lf}$ follows.
\\ \indent By Corollary \ref{convert} we know that
$\overline{C}^*$ and $\overline{C}$ are after a reflection polar
to each other. Therefore the have the same coefficient of symmetry
with respect to $r^{2k}$. Since $T_{2k,2k}$ fixes $r^{2k}$, it
follows that the coefficient of symmetry of $\overline{\lf}$ is
the same as the coefficient of symmetry of $\overline{C}^*$, which
by Theorem \ref{nonnegcoeffsym} is $(D(n,k)-1)^{-1}$.
\end{proof}
\begin{cor}
\label{powerball} $\overline{\lf}$ contains a ball of radius
\begin{equation*}
\frac{k!\Gamma(\frac{2k+n}{2})}
{\Gamma(\frac{4k+n}{2})\sqrt{D(n,k)-1}}
\end{equation*}
centered at $r^{2k}$.
\end{cor}
\begin{proof}
The coefficient
\begin{equation*}
\frac{(k-i)!\Gamma(\frac{2k+2i+n}{2})} {k!\Gamma(\frac{2k+n}{2})}
\quad \text{for} \quad 0 \leq i \leq k,
\end{equation*}
is clearly maximized when $i=k$. Thus the Loewner ellipsoid of
$\overline{\lf}$ contains a ball of radius
\begin{equation*}
\frac{k!\Gamma(\frac{2k+n}{2})\sqrt{D(n,2k)-1}}
{\Gamma(\frac{4k+n}{2})}.
\end{equation*}
From Proposition \ref{symmellip} and Theorem \ref{mainpowers} we
know that $\overline{\lf}$ will contain its Loewner ellipsoid
shrunk by the factor of
\begin{equation*}
\frac{1}{\sqrt{(D(n,2k)-1)(D(n,k)-1)}}.
\end{equation*}
Now the corollary follows.
\end{proof}
\begin{defin}
Let $K(2m)$ be the cone in $P_{n,2k}$ of forms whose restrictions
to the sphere are linear combinations of $2m$-th powers of linear
forms on $\sph$. Equivalently $K(2m)$ is the cone of forms in
$P_{n,2k}$ that multiplied by $r^{2m-2k}$ become sums of powers of
linear forms
\end{defin}
\begin{equation*}
K(2m)=\biggl{\{} f \in P_{n,2k} \mathrel{\bigg{\arrowvert}}
(x_1^2+\ldots+x_n^2)^{m-k}f \in \lf_{n,2m} \biggl{\}}.
\end{equation*}
We define $\overline{K}(2m)$ by intersecting $K$ with the
hyperplane of forms of integral 1 on $\sph$.
\begin{theorem}
\label{volumeratio} Let $\alpha=\text{dim}P_{n,2k}-1$. Then
\begin{equation*}
\left( \frac{\text{vol}\,\overline{K}(2m)}{\text{vol}\,
\overline{C}} \right)^{1/\alpha} \geq
\frac{m!\Gamma(\frac{2m+n}{2})} {(m-k)!\Gamma(\frac{2m+2k+n}{2})}.
\end{equation*}
\end{theorem}
\begin{proof}
We observe that from the definition of $T_{2m,2k}$ it follows that
$T_{2m,2k}$ maps $C$ into $K(2m)$. Since $T_{2m,2k}$ fixes
$r^{2k}$, it follows that $T_{2m,2k}$ maps $\overline{C}$ into
$\overline{K}(2m)$. But from Lemma \ref{rewrite} $T_{2m,2k}$ acts
on $H^*_{n,2i}$ by shrinking it by a factor of
\begin{equation*}
\frac{m!\Gamma(\frac{2m+n}{2})} {(m-i)!\Gamma(\frac{2m+2i+n}{2})}.
\end{equation*}
This coefficient is clearly minimized when $i=k$ and then the
theorem follows.
\end{proof}
\begin{cor}
Let $\epsilon \geq 0$ and let $m=(2k^2+kn)/{\epsilon}$. Then
\begin{equation*}
\left( \frac{\text{vol}\,\overline{K}(2m)}{\text{vol}\,
\overline{C}} \right)^{1/\alpha} \geq 1-\epsilon.
\end{equation*}
\end{cor}
\begin{remark}
The volume ratio
\begin{equation*}
\left( \frac{\text{vol}\,\overline{K}(2m)}{\text{vol}\,
\overline{C}} \right)^{1/\alpha}
\end{equation*}
approaches 1 as $m$ tends to infinity. Therefore every strictly
positive form lies in some $K(2m)$. (cf \cite{rez1})
\end{remark}

\section*{Acknowledgements}
The author wishes to thank Alexander Barvinok for many discussions
and support in writing this paper.


\begin{thebibliography}{dmst}

\bibitem[1]{ball} K. Ball,
\textit{An elementary introduction to modern convex geometry},
Flavors of Geometry 1-58, Math. Sci. Res. Inst. Publ., 31,
Cambridge Univ. Press, Cambridge, (1997).

\bibitem[2]{barv} A.I. Barvinok,
\textit{Estimating $L^{\infty}$ norms by $L^{2k}$ norms for
functions on orbits}, Foundations of Computational Mathematics, to
appear.

\bibitem[3]{realal}J. Bochnak, M. Coste, M-F. Roy,
\textit{Real Algebraic Geometry}, Springer-Verlag, Berlin, (1998).

\bibitem[4]{duo} J. Duoandkoetxea,
\textit{Reverse H\"{o}lder inequalities for spherical harmonics},
Proc. Amer. Math. Soc. 101 (1987), no. 3, 487-491.

\bibitem[4]{hilbert} D. Hilbert,
\textit{\"{U}ber die Darstellung definiter Formen als Summe von
Formenquadraten}, Math. Ann. 32, 342-350 (1888). Ges Abh. vol. 2,
415-436. Chelsea Publishing Co., New York, (1965)

\bibitem[5]{john} F. John,
\textit{Extremum problems with inequalities as subsidiary
conditions}, Studies and Essays Presented to R. Courant on his
60th Birthday, 187-204. Interscience Publishers, Inc., New York,
(1948).

\bibitem[6]{mull} C. M\"{u}ller,
\textit{Analysis of Spherical Symmetries in Euclidean Spaces},
Springer-Verlag, New York, (1998).

\bibitem[7]{rez1}B. Reznick,
\textit{Some concrete aspects of Hilbert's 17th Problem}, Contemp.
Math., 253 (2000), 251-272.

\bibitem[8]{rez2} B. Reznick,
\textit{Sums of even powers of real linear forms}, Mem. Amer.
Math. Soc. 96 (1992), no. 463.

\bibitem[9]{cgroup} B. Simon,
\textit{Representations of Finite and Compact Groups}, Graduate
Studies in Math, 10, AMS, Providence, RI (2001).

\bibitem[10]{sogge} C. Sogge,
\textit{Oscillatory integrals and spherical harmonics}, Duke Math.
J. 53 (1986), no. 1, 43-65.

\bibitem[11]{szeg} G. Szeg\"{o},
\textit{Orthogonal Polynomials}, American Mathematical Society
Colloquium Publications, v. 23., (1939).

\bibitem[12]{vilenkin} N. Ja. Vilenkin,
\textit{Special Functions and the Theory of Group Representations}
,Translations of Mathematical Monographs, Vol. 22, American
Mathematical Society (1968).


\end{thebibliography}
\end{document}